\documentstyle{amsppt}
\voffset-12mm
\magnification1200
\binoppenalty=10000\relpenalty=10000\relax
\let\hlop\!
\TagsOnRight
\addto\tenpoint{\normalbaselineskip=1.05\normalbaselineskip\normalbaselines}

\let\ge\geqslant

\let\wt\widetilde
\let\<\langle
\let\>\rangle
\let\kappa\varkappa
\define\trdeg{\operatorname{tr\,deg}}
\redefine\pmod#1{\;(\operatorname{mod}#1)}
\define\SL{\hbox{\sl SL}}
\define\GL{\hbox{\sl GL}}
\define\SP{\hbox{\sl Sp}}

\define\sL{\Cal L}
\define\sK{\Cal K}
\define\sC{\Cal C}
\define\fH{\frak H}
\redefine\d{\roman d}
\define\eprint#1#2{{\tt http://\allowlinebreak xxx.lanl.gov/\allowlinebreak
abs/\allowlinebreak#1/\allowlinebreak#2}}
\topmatter
\line{\hss math.NT/0008237}
\bigskip
\title
Number theory casting
a look at the mirror
\endtitle
\author
W.~Zudilin
\endauthor
\date
September 27, 2000
\enddate
\address
\hbox to70mm{\vbox{\hsize=70mm%
\leftline{Moscow Lomonosov State University}
\leftline{Department of Mechanics and Mathematics}
\leftline{Vorobiovy Gory, Moscow 119899 RUSSIA}
}}
\endaddress
\email
{\tt wadim\@ips.ras.ru}
\endemail
\dedicatory
To A.\,B.~Shidlovski\v\i\ on the occasion of his 85th birthday
\enddedicatory
\abstract
In this work, we give a purely analytic introduction
to the phenomenon of mirror symmetry for quintic threefolds
via classical hypergeometric functions
and differential equations for them.
Starting with a modular map and recent transcendence results
for its values, we regard a mirror map~$z(q)$ as a concept
generalizing the modular one. We give an alternative approach
demonstrating the existence of non-linear differential equations
for the mirror map, and exploit both an elegant construction
of Klemm--Lian--Roan--Yau and the Ax theorem to prove that the Yukawa
coupling~$K(q)$ does not satisfy any algebraic differential equation
of order less than~$7$ with coefficients from~$\Bbb C(q)$.
\endabstract
\endtopmatter
\document
\input mfpic
\opengraphsfile{ms}

It is a classical question of transcendence number theory
to investigate linear and algebraic independence of values
of analytic functions satisfying both arithmetic conditions
and functional (for instance, differential) equations.
This story has many dramatic and romantic episodes
(like the solution of the 7th Hilbert problem),
but it is far from an end.

Skipping results about the transcendence of one-variable modular
functions and their values, we concentrate
on what has been done in~1996.
The Mahler conjecture about the transcendence of at least one
among the numbers $q\in\Bbb C$, $0<|q|<1$, and $J(q)$, where
and $J(e^{2\pi i\tau})$~is a modular invariant,
has been proved in~\cite{BDGP}. Yu.~Nesterenko
has generalized this result using the Ramanujan functions
and differential equations for them; in~\cite{Ne} he has proved
that at least three among the numbers
$q$, $J(q)$, $\delta_qJ(q)$, and $\delta_q^2J(q)$,
for $q\in\Bbb C$ such that $0<|q|<1$ and $J(q)\notin\{0,1728\}$,
are algebraically independent over~$\Bbb Q$,
where $\delta_q=q\frac{\d}{\d q}$. It can be a subject
of an independent paper to overview consequences of the results
of~\cite{BDGP} and~\cite{Ne}; moreover, we believe that
the rich nature of modular functions will bring many new results
on the transcendence of mathematical constants.
The general theorem which allows one to apply the approach suggested
there to other analytic functions satisfying differential equations
has been also proved in~\cite{Ne}.
But all known (to us) non-linear differential systems
with `nice' (in some arithmetic sense) solutions
have the modular nature.

A recent observation of physicists in string theory
has given a job to mathematicians specializing in algebra
and geometry. This observation, called {\it mirror symmetry},
is nothing more than usual modularity in the simplest cases.
We are not planning here a lengthy discussion of the mirror referring
to wonderful articles on this subject (see~\cite{BS},
\cite{M1}--\cite{M3}).
Our aim is to work out an approach to mirror symmetry as a philosophy
generalizing modular, in particular, as a good task
for number theorists and analysts. Below we make an attempt
to describe mirror symmetry analytically, without any physics
and algebraic geometry.

\head
1. Modularity as a motivation
\endhead

Let us start from the solution~$f_0(z)$ analytic at $z=0$
of the second order linear differential equation
$$
\bigl(\delta_z^2-3z(3\delta_z+1)(3\delta_z+2)\bigr)y=0,
\qquad\text{where}\quad \delta_z=z\frac{\d}{\d z}.
\tag1
$$
One finds that $f_0(z)$ is the hypergeometric $G$-function
with expansion
$$
f_0(z)
=\sum_{l=0}^\infty\frac{(3l)!}{(l!)^3}z^l
={}_2\!F_1\biggl(\frac13,\frac23;1;3^3z\biggr)
$$
in a neighbourhood of $z=0$.
Since the equation~\thetag{1} has unipotent monodromy
at the point $z=0$, another solution of~\thetag{1}
is of the form $f_1(z)=f_0(z)\log z+g(z)$, where
$$
g(z)=3\sum_{l=1}^\infty
\biggl(\frac{(3l)!}{(l!)^3}\sum_{k=l+1}^{3l}\frac1k\biggr)z^l
$$
in a neighbourhood of $z=0$.
It is non-obvious to see that the function
$$
\align
q(z)&=\exp\biggl(\frac{f_1(z)}{f_0(z)}\biggr)
=z\cdot\exp\biggl(\frac{g(z)}{f_0(z)}\biggr)
\\
&=z+15z^2+279z^3+5729z^4+124554z^5+2810718z^6+65114402z^7+O(z^8)
\endalign
$$
has integral coefficients in its $z$-expansion,
so the same property holds for the $q$-expansion of the inverse function
$$
z(q)=q-15q^2+171q^3-1679q^4+15054q^5-126981q^6+1024952q^7+O(q^8).
$$
Moreover, we have the expansion
$$
\wt f_0(q)=f_0\bigl(z(q)\bigr)
=1+6q+6q^3+6q^4+12q^7+6q^9+6q^{12}+12q^{13}+6q^{16}+12q^{19}+O(q^{21}),
$$
which converges for $q\in\Bbb C$ with $|q|<1$, and
$$
f_0^2\bigl(z(q)\bigr)
=\biggl(\frac{\delta_qz(q)}{z(q)}\biggr)\cdot\frac1{1-3^3z(q)}
\tag2
$$
(see, e.g., \cite{LY2} or~\cite{Z}).
These observations can be explained by the modular origin
of functions $z(q)$ and $\wt f_0(q)$ with respect to the parameter
$\tau=\frac1{2\pi i}\log q\in\fH
=\{\tau\in\Bbb C:\allowmathbreak\Im\tau>0\}$
and the congruence level~$3$ subgroup
$$
\Gamma_0(3)=\biggl\{\pmatrix a&b\\c&d\endpmatrix\in\SL_2(\Bbb Z):
c\equiv0\pmod3\biggr\}
$$
(see~\cite{LY2},~\cite{HM},~\cite{Z}): each of these functions satisfies
a third order non-linear algebraic equation over~$\Bbb Q$
with respect to $\delta_q$-derivation,
where $\delta_q=q\frac{\d}{\d q}=\frac1{2\pi i}\frac{\d}{\d\tau}$,
and, due to Mahler~\cite{Ma} (see also~\cite{Ni}),
does not satisfy any equation over~$\Bbb C(q)$ of a smaller order.

Note that $f_0(z)$ and $q(z)=\exp(f_1(z)/f_0(z))$
are multi-valued functions of
$z\in\allowmathbreak\Bbb C\setminus\{3^{-3}\}$;
the function~$z(q)$ is analytic for $q\in\Bbb C$ with $|q|<1$
except simple poles at points~$q$ such that
$\tau=\frac1{2\pi i}\log q$ is congruent to
$\frac1{\sqrt3}e^{\pi i/6}=\frac12+\frac i{2\sqrt3}$
with respect to~$\Gamma_0(3)$
(see the fundamental domain of this group on Fig.~1).
But locally, in a neighbourhood of $z\in\Bbb C\setminus\{0,3^{-3}\}$,
for any fixed branches of $f_0(z)$ and $q(z)$, the fields
$$
\Bbb Q\bigl(z,q(z),f_0(z),\delta_zf_0(z)\bigr)
\quad\text{and}\quad
\Bbb Q\bigl(q,z(q),\wt f_0(q),\delta_q\wt f_0(q)\bigr)
\tag3
$$
coincide up to algebraic extension.

\midinsert
\line{\hss
\mfpic[15]{-6.3}{6.3}{-.5}{10}
\pointsize3pt \dashlen3pt \dashspace1.5pt
\headshape{1}{4}{false}
\shade[18]\closed\connect
\lines{(3,10),(3,1.732)}
\arc[p]{(2,0),60,180,2}
\arc[p]{(-2,0),0,120,2}
\lines{(-3,1.732),(-3,10)}
\endconnect
\point{(3,1.732),(0,0),(-3,1.732)}
\label[cl](3.3,2.2){$\frac12+\frac i{2\sqrt3}$}
\label[tc](0,-.25){$0$}
\label[cr](-3.3,2.2){$-\frac12+\frac i{2\sqrt3}$}
\white\ellipse{(.5,9.5),.8,.54}
\arrow\lines{(0,0),(0,10)}
\label[cc](.6,9.5){$\Im\tau$}
\arrow\lines{(-6.3,0),(6.3,0)}
\label[cc](5.7,.5){$\Re\tau$}
\dotted\arc[p]{(2,0),0,60,2}
\dotted\arc[p]{(-2,0),120,180,2}
\pen{1pt}
\lines{(3,10),(3,1.732)}
\arc[p]{(2,0),60,180,2}
\arc[p]{(-2,0),0,120,2}
\lines{(-3,1.732),(-3,10)}
\endmfpic
\hss}
\botcaption{Fig.~1}
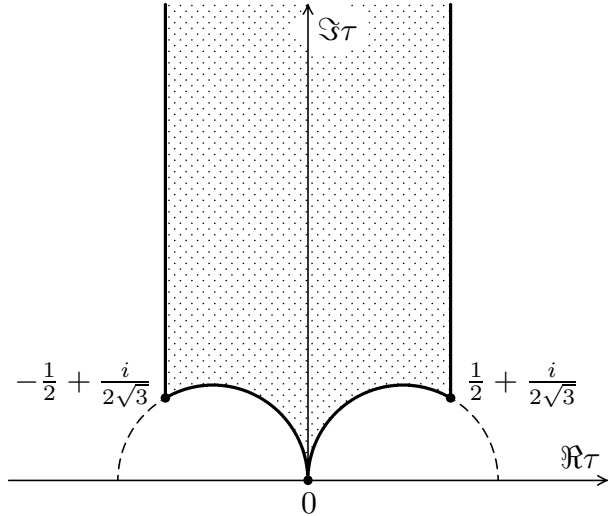
Fundamental domain of $\Gamma_0(3)$
\endcaption
\endinsert

One consequence of Nesterenko's result~\cite{Ne} cited above
is the algebraic independence of at least three among the numbers
$q$, $z(q)$, $\wt f_0(q)$, and $\delta_q\wt f_0(q)$
for each~$q$, $0<|q|<1$, such that $z(q)\notin\{0,3^{-3}\}$.
From the coincidence (up to algebraic extension) of the fields~\thetag{3}
we deduce the following result.

\proclaim{Theorem 1}
For each $z\in\Bbb C\setminus\{0,3^{-3}\}$\rom,
at least three numbers among
$z$\rom, $q(z)$\rom, $f_0(z)$\rom, and $\delta_zf_0(z)$
are algebraically independent over~$\Bbb Q$.
\endproclaim

Note that there is no general approach to
transcendence (and even irrationality!) proofs for values
of $G$-functions.

Now, we go on to the third order linear differential equation
$$
\bigl(\delta_z^3-4z(4\delta_z+1)(4\delta_z+2)(4\delta_z+3)\bigr)y=0
\tag4
$$
and take two its solutions $f_0(z)$ analytic at $z=0$
and $f_1(z)=f_0(z)\log z+g(z)$, where
$$
\gather
f_0(z)
=\sum_{l=0}^\infty\frac{(4l)!}{(l!)^4}z^l
={}_3\!F_2\biggl(\frac14,\frac24,\frac34;1,1;4^4z\biggr),
\\
g(z)=4\sum_{l=1}^\infty
\biggl(\frac{(4l)!}{(l!)^4}\sum_{k=l+1}^{4l}\frac1k\biggr)z^l
\endgather
$$
in a neighbourhood of $z=0$. Once again, the function
$$
\align
q(z)&=\exp\biggl(\frac{f_1(z)}{f_0(z)}\biggr)
=z\cdot\exp\biggl(\frac{g(z)}{f_0(z)}\biggr)
\\
&=z+104z^2+15188z^3+2585184z^4+480222434z^5+94395247376z^6+O(z^7)
\endalign
$$
has integral coefficients in its $z$-expansion,
and the inverse function
$$
z(q)=q-104q^2+6444q^3-311744q^4+13018830q^5-493025760q^6+O(q^7)
$$
has integral coefficients in its $q$-expansion;
moreover, the $q$-expansion
$$
\wt f_0(q)=f_0\bigl(z(q)\bigr)
=1+24q+24q^2+96q^3+24q^4+144q^5+96q^6+192q^7+24q^8+O(q^9)
$$
with integral coefficients converges for $q\in\Bbb C$, $|q|<1$,
and
$$
f_0^2\bigl(z(q)\bigr)
=\biggl(\frac{\delta_qz(q)}{z(q)}\biggr)^2\cdot\frac1{1-4^4z(q)}
\tag5
$$
(see~\cite{LY2}).
Theorem~1 (with the change~$3^{-3}$ by~$4^{-4}$)
holds for these new functions, since
$$
f_0(z)
=\sum_{l=0}^\infty\frac{(4l)!}{(l!)^4}z^l
=\biggl({}_2\!F_1\biggl(\frac18,\frac38;1;4^4z\biggr)\biggr)^2
$$
(see~\cite{LY2},~\cite{Z})
and equation~\thetag{4} is the symmetric square of
the second order linear differential equation
$$
\bigl(\delta_z^2-4z(8\delta_z+1)(8\delta_z+3)\bigr)y=0,
$$
which has a modular explanation like the equation~\thetag{1}.

\head
2. Natural generalization
\endhead

It is natural to continue the modular story of the previous section
by considering now the linear differential equation
$$
\bigl(\delta_z^{s-1}
-sz(s\delta_z+1)(s\delta_z+2)\dotsb(s\delta_z+s-1)\bigr)y=0,
\tag6
$$
where $s\ge3$ is an integer. Writing
$$
f(z;H)=z^H\sum_{l=0}^\infty z^l
\frac{\prod_{k=1}^{sl}(sH+k)}{\prod_{k=1}^l(H+k)^s}\pmod{H^{s-1}},
\tag7
$$
where
$$
z^H=e^{H\log z}=1+H\log z+H^2\frac{\log^2z}2+\dots
+H^{s-2}\frac{\log^{s-2}z}{(s-2)!}\pmod{H^{s-1}},
$$
one can verify that the functions $f_0,f_1,\dots,f_{s-2}$ from
the formal expansion
$$
f(z;H)=f_0(z)+f_1(z)H+\dots+f_{s-2}(z)H^{s-2}
$$
form the fundamental solution to the equation~\thetag{6}.
By~\thetag{7}, we have $z$-expansions
$$
\align
f_0(z)&=g_0(z)
=\sum_{l=0}^\infty\frac{(sl)!}{(l!)^s}z^l,
\tag8
\\
f_1(z)&=g_0(z)\log z+g_1(z),
\qquad\text{where}\quad
g_1(z)=s\sum_{l=1}^\infty
\biggl(\frac{(sl)!}{(l!)^s}\sum_{k=l+1}^{sl}\frac1k\biggr)z^l
\endalign
$$
in a neighbourhood of $z=0$.

The integrality of the expansions for
$$
q(z)=\exp\biggl(\frac{f_1(z)}{f_0(z)}\biggr)
=z\cdot\exp\biggl(\frac{g_1(z)}{g_0(z)}\biggr)
$$
and, as a consequence, for the inverse function~$z(q)$
has been recently proved by Lian and Yau~\cite{LY2}
using Dwork's $p$-adic approach. We underline that
no modular argument is known for~$z(q)$ when $s\ge5$;
moreover, for $s\ge5$, it is easy to show that $z(e^{2\pi i\tau})$ cannot be
a modular function with respect to a congruence subgroup of~$\SL_2(\Bbb Z)$
(that is a consequence of the structure of monodromy
groups for equations like~\thetag{6}; see~\cite{BH}
and Section~5 below; see also~\cite{D} for the answer
to the question `When is~$z(q)$ a modular function?').
We call $q\mapsto z(q)$ the {\it mirror map\/} produced
by differential equation~\thetag{6}, or by hypergeometric
series~\thetag{8}. (One can see no
real mirror in the construction above, but we would like
to have a compact name for the map~$z(q)$.)

Although the algorithm for deducing differential equations for
a mirror map is known (see~\cite{LY1}), we give another approach
via the {\it Wronskian formalism\/} in the next section.
In Section~4, we make a special emphasis on the extremely well studied
non-modular case of $s=5$. Finally, in Section~5, we prove
some transcendence results for the mirror map and related functions.

Throughout this paper we use the following convention
for the parameters $t,\tau,q$:
$$
t=\log q=2\pi i\tau,
$$
and for the corresponding derivations:
$$
\delta_q=q\frac{\d}{\d q}=\frac{\d}{\d t}
=\frac1{2\pi i}\frac{\d}{\d\tau}.
$$

\head
3. Wronskian formalism
\endhead

In this section, we derive the algebraic differential equation
satisfied by the inverse of a ratio of two independent
solutions of a linear differential equation (with regular singular points).
In the case of second order differential equations
such is the {\it Schwarzian equation\/}
$$
2Q(z)\biggl(\frac{\d z}{\d t}\biggr)^2+\{z,t\}=0,
\tag9
$$
where
$$
\{z,t\}=\frac{\d^3z/\d t^3}{\d z/\d t}
-\frac32\biggl(\frac{\d^2z/\d t^2}{\d z/\d t}\biggr)^2
$$
is the Schwarzian derivative, and
$Q(z)$ is a rational function
with poles of order at most two at the singular points.
Note that scaling the parameter~$t$ in~\thetag{9}
has no effect on this equation; so, it remains true
if we change $t$ by~$\tau$.

We are following the approach of Chudnovsky~\cite{CC}.
To deduce the desired equation for a general
$m$th order linear differential equation
$$
D[y]
=\frac{\d^my}{\d z^m}+a_1(z)\frac{\d^{m-1}y}{\d z^{m-1}}+\dots
+a_{m-1}(z)\frac{\d y}{\d z}+a_m(z)y
=0,
\tag10
$$
where the coefficients $a_1,\dots,a_m$ are from
a differentially closed field~$\sL$ (usually, $\sL=\Bbb C(z)$)
with constant field~$\Bbb C$,
we need some properties of the Wronskian
$$
W(f_0,\dots,f_{m-1})
=\det\biggl(\frac{\d^kf_j}{\d z^k}\biggr)_{k,j=0,1,\dots,m-1}.
$$
By~$\sL_1$ we denote a differential extension of~$\sL$
with the same constant field~$\Bbb C$.

\proclaim{Lemma 1}
The Wronskian $W(f_0,\dots,f_{m-1})$ of elements
$f_0,\dots,f_{m-1}\in\sL_1\supset\sL$
is identically zero if and only if
these elements are linearly dependent over the constant field~$\Bbb C$.
\endproclaim

\proclaim{Lemma 2}
Let $f_0,\dots,f_{m-1}\in\sL_1\supset\sL$
be linearly independent over~$\Bbb C$ solutions of~\thetag{10}.
Then
$$
D[y]=\frac{W(y,f_0,\dots,f_{m-1})}{W(f_0,\dots,f_{m-1})}.
$$
\endproclaim

\proclaim{Lemma 3}
For any collection $g,f_0,\dots,f_{m-1}\in\sL_1$\rom,
$$
W(gf_0,\dots,gf_{m-1})=g^mW(f_0,\dots,f_{m-1}).
$$
\endproclaim

First, for given $m$th order linear differential equation~\thetag{10}
we want to construct the (non-linear) differential equation
satisfied by each ratio of two linearly independent solutions
of~\thetag{10}. It is
$$
R[t]
=\frac{W(tf_0,tf_1,\dots,tf_{m-1},f_0,f_1,\dots,f_{m-1})}
{W^2(f_0,f_1,\dots,f_{m-1})}
=0,
$$
where $f_0,f_1,\dots,f_m$ are linearly independent solutions
of~\thetag{10}. By Lemmas~1 and~2, the coefficients of $R[t]$
belong to the differentially closed field~$\sL$ containing
all coefficients of the equation~\thetag{10};
the order of differential operator~$R[t]$ is~$2m-1$.

Using the well-known binomial formula for the multiple derivation
of the product and eliminating higher-order derivatives of
$f_0,f_1,\dots,f_{m-1}$ from the last $m$~columns
in the numerator of~$R[t]$, we get
$$
R[t]
=\frac1{W(f_0,\dots,f_{m-1})}
\det\biggl(\sum_{l=1}^k\binom klt^{(l)}
\frac{\d^{k-l}f_j}{\d z^{k-l}}\biggr)_{k=m,m+1,\dots,2m-1;\ j=0,1,\dots,m-1},
$$
which yields that $R[t]$ is a polynomial in $t',\dots,t^{(2m-1)}$
with coefficients from~$\sL$ such that its leading term $t^{(2m-1)}$
has coefficient~$1$.

In the case $\sL=\Bbb C(z)$ we make a transition from the variable~$z$
and the function~$t(z)$ to the variable~$t$ and the function~$z(t)$.
Since the fraction fields
$$
\Bbb Q\biggl(\frac{\d t}{\d z},\frac{\d^2t}{\d z^2},\dots,
\frac{\d^{2m-1}t}{\d z^{2m-1}}\biggr)
\quad\text{and}\quad
\Bbb Q\biggl(\frac{\d z}{\d t},\frac{\d^2z}{\d t^2},\dots,
\frac{\d^{2m-1}z}{\d t^{2m-1}}\biggr)
$$
coincide, $R[t]$ becomes a rational
function of $z',z'',\dots,z^{(2m-1)}$ with coefficients from~$\Bbb C(z)$
or, equivalently, $R[t]=0$ becomes a polynomial differential
equation of $(2m-1)$th order with constant coefficients.
Moreover, $z^{(2m-1)}$ enters this equation in the first power.

Applying this construction to the case of the linear differential
equation~\thetag{6} we obtain the following result.

\proclaim{Proposition 1}
The mirror map produced by differential equation~\thetag{6}
satisfies an algebraic $(2s-3)$th order differential equation
with coefficients from~$\Bbb C$. Moreover\rom, $z^{(2s-3)}$~is
rational over the field $\Bbb C(z,z',\dots,z^{(2s-4)})$.
\endproclaim

\remark{Remark}
The algebraic equation constructed above coincides with
the equation given by the algorithm from~\cite{LY1}.
\endremark

Now we present some preliminary results for our transcendence consideration
in Section~5. Let $f_0,f_1,\dots,f_{m-1}$~be
a fundamental solution of the (Fuchsian) differential
equation~\thetag{10} with coefficients from~$\Bbb C(z)$,
in a neighbourhood of a non-singular point $z=z_0$.
Then each of the functions
$$
t_j(z)=\frac{f_j(z)}{f_0(z)},
\qquad j=1,\dots,m-1,
\tag11
$$
satisfies the same algebraic differential equation $R[t]=0$.
If $\sL$~denotes the Picard--Vessiot extension of~$\Bbb C(z)$
corresponding to~\thetag{10}, then
from the definition~\thetag{11} we see that $t_j\in\sL$
for all $j=1,\dots,m-1$. Further, the field
$$
\sK=\Bbb C\biggl(z,t_j,\frac{\d t_j}{\d z},\dots,
\frac{\d^{2m-2}t_j}{\d z^{2m-2}}\biggr)_{j=1,\dots,m-1}
\subset\sL
\tag12
$$
is differentially stable. By Lemma~3,
$$
f_0^m=\frac{W(f_0,f_1,\dots,f_{m-1})}{W(1,t_1,\dots,t_{m-1})},
\tag13
$$
and we note that the numerator of~\thetag{13}---the
Wronskian of~\thetag{10}---is an algebraic function since
the equation~\thetag{10} is Fuchsian. Thus, by~\thetag{13}
we derive that $f_0$~is algebraic over~$\sK$; hence the same
property holds for $f_j=t_jf_0$, $j=1,\dots,m-1$.
Therefore we have

\proclaim{Proposition 2}
The differentially stable fields $\sK$ and $\sL$ defined above
coincide up to algebraic extension.
\endproclaim

\remark{Remark}
As an easy consequence of Lemma~1
it is possible to deduce the {\it linear\/} equation satisfied
by the functions~\thetag{11}, but the coefficients of this equation
are not in~$\Bbb C(z)$.
\endremark

\head
4. Differential equations for mirror map
\endhead

In this section we make a special emphasis on the case $s=5$
in~\thetag{6}. We start with a comparison of expansions with those
holding in the modular cases of Section~1.
In the case $s=5$ we derive the $q$-expansions
$$
\align
z(q)
&=q-770q^2+171525q^3-81623000q^4-35423171250q^5-54572818340154q^6
\\ &\qquad
-71982448083391590q^7-102693620674349200800q^8+O(q^9),
\tag14
\\
\wt f_0(q)
&=f_0\bigl(z(q)\bigr)
\\
&=1+120q+21000q^2+14115000q^3+13414125000q^4+15234972675120q^5
\\ &\qquad
+19285869813670920q^6+26264963911492602000q^7+O(q^8).
\tag15
\endalign
$$
Unfortunately, the convergence domain of the $q$-expansion~\thetag{15}
is not the disc $|q|<1$.

Without using the results of Section~3, we now present a compact version
of the $\delta_q$-differential equation for the function $z(q)$
(\cite{LY1}--\cite{LY3}):
$$
2Q(z)\biggl(\frac{\d z}{\d t}\biggr)^2+\{z,t\}
=\frac25\frac{\d^2\log K}{\d t^2}
-\frac1{10}\biggl(\frac{\d\log K}{\d t}\biggr)^2
\tag16
$$
(cf.~\thetag{9}), where
$$
\align
Q(z)
&=\frac{5^8}4\biggl(\frac{16}{5^5z}+\frac{16}{1-5^5z}
+\frac{25}{(5^5z)^2}+\frac{15}{(1-5^5z)^2}\biggr)
\\
&=\frac{5^8}4\cdot\frac{25-34\cdot(5^5z)+24\cdot(5^5z)^2}
{(5^5z)^2(1-5^5z)^2}
\tag17
\endalign
$$
is a rational function and
$$
\align
K(q)
&=5+\sum_{l=1}^\infty\frac{n_ll^3q^l}{1-q^l}
\\
&=5+2875q+4876875q^2+8564575000q^3+15517926796875q^4+O(q^5)
\tag18
\endalign
$$
is the so-called {\it Yukawa coupling}.

We cannot write the direct analogue of~\thetag{2} and~\thetag{5}
in this case, but we have
$$
f_0^2\bigl(z(q)\bigr)
=\biggl(\frac{\delta_qz(q)}{z(q)}\biggr)^3\cdot\frac1{1-5^5z(q)}
\cdot\frac5{K(q)}
\tag19
$$
(see, e.g., \cite{BS}). Formula~\thetag{19} can be regarded
as the definition of the Yukawa coupling, although the original definition
is based on the enumerative meaning of the numbers~$n_l$ in~\thetag{18}
(see~\cite{Pa}). The definition~\thetag{19} and the integrality
of the mirror map~\thetag{14} yield the conclusion
$\frac15K(q)\in\Bbb Z[[q]]$
without any references to algebraic geometry.

Picking, as in~\thetag{7}, the fundamental solution
$f_0,f_1,f_2,f_3$ of the equation
$$
\bigl(\delta_z^4
-5z(5\delta_z+1)(5\delta_z+2)(5\delta_z+3)(5\delta_z+4)\bigr)y=0
\tag20
$$
we get
$$
\aligned
f_0(z)&=g_0(z),
\qquad
f_1(z)=g_0(z)\log z+g_1(z),
\\
f_2(z)&=g_0(z)\frac{\log^2z}2+g_1(z)\log z+g_2(z),
\\
f_3(z)&=g_0(z)\frac{\log^3z}6+g_1(z)\frac{\log^2z}2+g_2(z)\log z+g_3(z),
\endaligned
$$
where
$$
\align
g_0(z)&=1+120z+113400z^2+168168000z^3+305540235000z^4+O(z^5),
\\
g_1(z)&=770z+810225z^2+\frac{3745679000}3z^3+\frac{4627120640625}2z^4+O(z^5),
\\
g_2(z)&=575z+\frac{4208175}4z^2+\frac{16964522000}9z^3
+\frac{180021646778125}{48}z^4+O(z^5),
\\
g_3(z)&=-1150z-\frac{3298375}4z^2-\frac{46661619875}{54}z^3
-\frac{325329574909375}{288}z^4+O(z^5)
\endalign
$$
are analytic functions. Taking
$$
t(z)=\frac{f_1(z)}{f_0(z)}
=\log q(z)=2\pi i\tau(z)
\tag21
$$
for a new variable and
$$
t_j(t)=\frac{f_j(z(t))}{f_0(z(t))}, \qquad j=0,1,2,3,
\tag22
$$
for the new functions,
we obtain for~\thetag{22} the differential equation
$$
\frac{\d^2}{\d t^2}\frac1{K(e^t)}\frac{\d^2}{\d t^2}t_j=0,
\qquad j=0,1,2,3
$$
(see~\cite{Pa}), where $K$~is the Yukawa coupling. Moreover,
we can explicitly describe functions~\thetag{22}.
Namely, taking a {\it prepotential\/}
$$
F(t)=\frac56t^3+\sum_{l=1}^\infty N_lq^l
=\frac56t^3+\sum_{l=1}^\infty\sum_{k=1}^\infty\frac{n_lq^{kl}}{k^3},
$$
and using the argument of~\cite{LY3}
we obtain for the equation~\thetag{20}
$$
t_0(t)=1, \quad t_1(t)=t \quad \text{(that is obvious)}, \qquad
t_2(t)=\frac15\frac{\d F}{\d t}, \quad
t_3(t)=\frac15t\frac{\d F}{\d t}-\frac25F,
\tag23
$$
which, combined with~\thetag{21} and~\thetag{22}, yields
$$
F\biggl(\frac{f_1(z)}{f_0(z)}\biggr)
=\frac52\frac{f_1(z)f_2(z)-f_0(z)f_3(z)}{f_0^2(z)}.
$$
Moreover, it is easy to see that
$$
\frac{\d^3F}{\d t^3}(t)=K(e^t)=K(q).
\tag24
$$

\remark{Remark}
The functions
$$
F_0(t)=\frac16t^3
+\sum_{l=1}^\infty\sum_{k=1}^\infty\frac{240q^{kl}}{k^3},
\qquad
K_0(q)=\frac{\d^3F_0}{\d t^3}(t)
=1+\sum_{l=1}^\infty\frac{240l^3q^l}{1-q^l}
$$
can be viewed as the `modular' analogues of $F(t)$ and~$K(q)$.
The function~$K_0(q)$~is a well known Ramanujan function,
or an Eisenstein series regarded as a function of
$\tau=\allowmathbreak\frac1{2\pi i}\log q$.
Among properties of~$F_0(t)$ we mention an arithmetic one:
for $t=2\pi i\tau=\allowmathbreak-2\pi$,
$$
F_0(-2\pi)=\frac{10}3\pi^3-120\zeta(3),
\qquad \text{where} \quad
\zeta(3)=\sum_{l=1}^\infty\frac1{l^3}
$$
(see footnoted Ramanujan's identity in~\cite{Po}).
\endremark

Now, we want to expose the results from~\cite{KLRY} and~\cite{LY3}
on the existence of a differential equation for the Yukawa coupling.
Moreover, as shown in these papers, there exists
a duality between the differential equations for the mirror map
and for the Yukawa coupling.

Any fourth order linear differential equation producing
a mirror map with integral expansion can be reduced to the form
$$
\frac{\d^4y}{\d z^4}+Q_2(z)\frac{\d^2y}{\d z^2}
+\frac{Q_2(z)}{\d z}\frac{\d y}{\d z}+Q_0(z)y=0
$$
(see~\cite{LY3}), where $Q_2(z)$ and $Q_0(z)$ are rational functions
uniquely determined by the original equation (for instance,
in the case of equation~\thetag{20}
one has $Q_2(z)=\allowmathbreak 10Q(z)$
with~$Q$ defined in~\thetag{17}). Using primes for
$\d/\d t$-derivatives and following \cite{KLRY},~\cite{LY3} we set
$$
\aligned
A_2(z)&=Q_2(z){z'}^2+5\{z;t\},
\\
A_4(z)&=Q_0(z){z'}^4+\frac32\frac{\d Q_2(z)}{\d z}{z'}^2z''
-\frac34Q_2(z){z''}^2
+\frac32Q_2(z)z'z^{(3)}
\\ &\qquad
-\frac{135}{64}\biggl(\frac{z''}{z'}\biggr)^4
+\frac{75}4\frac{{z''}^2z^{(3)}}{{z'}^3}
-\frac{15}4\biggl(\frac{z^{(3)}}{z'}\biggr)^2
-\frac{15}2\frac{z''z^{(4)}}{{z'}^2}+\frac32\frac{z^{(5)}}{z'},
\\
B_2(u)&=2u''-\frac{{u'}^2}2,
\qquad
B_4(u)=\frac{u^{(4)}}2+\frac{{u''}^2}4-\frac{u''{u'}^2}2+\frac{{u'}^4}{16}.
\endaligned
$$

\widestnumber\item{\kern1mm}
\proclaim{Proposition 3 \cite{KLRY}}
Given a pair of rational functions $Q_0(z),Q_2(z)$
\rom(which determines the Picard--Fuchs equation\rom)\rom,
there exists a differential polynomial $P_1(\,\cdot\,,\,\cdot\,)$
with the following properties\rom:
\roster
\item"(i)" $P_1$ is quasi-homogeneous\rom;
\item"(ii)" $P_1(A_2(z),A_4(z))$ is identically zero\rom;
\item"(iii)" $P_1(B_2(u),B_4(u))=0$ is a non-trivial seventh order
differential equation in~$u$ with a solution $u(t)=\log K(e^t)$\rom;
\item"(iv)" $P_1$ is minimal\rom, i.e\.
any differential polynomial satisfying~\rom{(i)--(iii)}
has degree not less than~$P$\rom;
\item"(v)" the differential equation $P_1(B_2(u),B_4(u))=0$
is $\SL_2(\Bbb C)$-invariant \rom(that is\rom, it is stable
under a change of variable $t\mapsto(at+b)/(ct+d)$ for
$\bigl(\smallmatrix a&b\\c&d\endsmallmatrix\bigr)\in\SL_2(\Bbb C)$\,\rom).
\endroster
The polynomial $P_1$ characterized by the properties \rom{(i)--(v)}
is unique up to a constant multiple.
\endproclaim

We do not present here the complicated differential equation
for the Yukawa coupling \thetag{18} described in Proposition~3,
referring to an elegant algorithm for its construction and another example
in~\cite{KLRY} and~\cite{LY3}.

\widestnumber\item{\kern1mm}
\proclaim{Proposition 4 \cite{KLRY}}
There exists a differential polynomial $P_2(\,\cdot\,,\,\cdot\,)$
with the following properties\rom:
\roster
\item"(i)" $P_2$ is quasi-homogeneous\rom;
\item"(ii)" $P_2(B_2(u),B_4(u))$ is identically zero\rom;
\item"(iii)" $P_2(A_2(z),A_4(z))=0$ is a non-trivial seventh order
differential equation in~$z$ with a solution $z(e^t)$\rom;
\item"(iv)" $P_2$ is minimal of degree~$12$\rom, i.e\.
any differential polynomial satisfying~\rom{(i)--(iii)}
has degree at least~$12$\rom;
\item"(v)" the differential equation $P_2(A_2(z),A_4(z))=0$
is $\SL_2(\Bbb C)$-invariant\rom;
\item"(vi)" $P_2$ is universal\rom, i.e\.
it is independent of the data $Q_0(z),Q_2(z)$
and it is characterized by the properties \rom{(i)--(v)}
up to a constant multiple.
\endroster
The equation $P_2(A_2(z),A_4(z))=0$ coincides with the equation
deduced in Sec\-tion~\rom3.
\endproclaim

As noted in~\cite{LY1} and proved in~\cite{LY3},
it is possible to deduce coupled non-linear differential equations
solved by the mirror map and the Yukawa coupling, e.g.,
\thetag{16}~and
$$
\wt Q(z)\biggl(\frac{z'}z\biggr)^4
=\frac{175{K'}^4-280K{K'}^2K''+49K^2{K''}^2
+70K^2K'K^{(3)}-10K^3K^{(4)}}{K^4},
\tag25
$$
where $\wt Q(z)$ is a rational function depending on
the original linear differential equation. In the case~\thetag{20},
$$
\align
\wt Q(z)
&=-\frac{5750z+63671875z^2+19531250000z^3}{(1-3125z)^4}
\\
&=-\frac1{25}\frac{2\cdot5^5z+163\cdot(5^5z)^2+8\cdot(5^5z)^3}{(1-5^5z)^4}.
\endalign
$$

\proclaim{Proposition 5}
Let $\sK_0$ be the algebraic closure over~$\Bbb C$ of the field
generated by the function $K(q)$ and its $\delta_q$-derivatives,
where $\delta_q=q\frac{\d}{\d q}=\frac{\d}{\d t}$.
Then the mirror map $z(q)$ is algebraic over~$\sK_0$.
\endproclaim

\demo{Proof}
By Proposition~3, the transcendence degree of~$\sK_0$ over~$\Bbb C$
is at most~$7$. We extend algebraically the field~$\sK_0$ to~$\wt{\sK}_0$
by adding the element~$\mu$ that is the root of the fourth degree of the
right-hand side of~\thetag{25}; then the derivatives $\mu',\mu'',\dots$
also belong to~$\wt{\sK}_0$.
It is sufficient to prove that $z(q)$ is algebraic over~$\wt{\sK}_0$.
To this end we rewrite \thetag{25} as follows:
$$
z'=R(z)\cdot\mu, \qquad \mu\in\wt{\sK}_0,
\tag26
$$
where $R(z)$ is an algebraic function of~$z$, and take the logarithmic
derivative of~\thetag{26}:
$$
\frac{z''}{z'}=\frac{z'R'(z)}{R(z)}+\frac{\mu'}{\mu}
=R'(z)\mu+\mu_1,
\tag27
$$
where we set $\mu_1=\mu'/\mu\in\wt{\sK}_0$ and $R'(z)$ means
derivative of~$R(z)$ with respect to its intrinsic parameter~$z$.
Further, by~\thetag{27} we obtain
$$
\align
\{z;t\}
&=\biggl(\frac{z''}{z'}\biggr)'-\frac12\biggl(\frac{z''}{z'}\biggr)^2
\\
&=\bigl(z'R''(z)\mu+R'(z)\mu\mu_1+\mu_1'\bigr)
-\frac12\bigl({R'}^2(z)\mu^2+2R'(z)\mu\mu_1+\mu_1^2\bigr)
\\
&=\mu^2\Bigl(R(z)R''(z)-\frac12{R'}^2(z)\Bigr)
+\Bigl(\mu_1'-\frac12\mu_1^2\Bigr).
\tag28
\endalign
$$
On the other hand, by~\thetag{16} we obtain
$$
\{z;t\}
=-2Q(z){z'}^2+\mu_2
=-2Q(z)R^2(z)\mu^2+\mu_2,
\tag29
$$
where $\mu_2\in\sK_0\subset\wt{\sK}_0$ is the right-hand side of~\thetag{16}.
Comparing \thetag{28} and \thetag{29} we see that
$$
R(z)R''(z)-\frac12{R'}^2(z)+2Q(z)R^2(z)
=\frac1{\mu^2}\biggl(\frac12\mu_1^2-\mu_1'+\mu_2\biggr)
\in\wt{\sK}_0.
\tag30
$$
Direct calculations show that the left-hand side of~\thetag{30}
is not identically zero, which means the algebraicity of~$z$ over~$\wt{\sK}_0$
and, finally, over~$\sK_0$. This completes the proof.
\enddemo

\head
5. Transcendence problems of mirror map
\endhead

We are now able to state results on functional transcendence
for the mirror map and for the Yukawa coupling. To do this,
we apply the general method introduced in the joint works~\cite{BZ1}
and \cite{BZ2} of D.~Bertrand and this author.
We go back to notation of Section~3, where $\sL$~is the Picard--Vessiot
extension corresponding to equation~\thetag{20}
and the field~$\sK$ from~\thetag{12} is $\d/\d z$-differentially stable.

Let $t=t_1$ be a new parameter. Since
$$
\frac{\d}{\d t}=\biggl(\frac{\d t}{\d z}\biggr)^{-1}\frac{\d}{\d z},
$$
we can consider the field~$\sK$ as $\d/\d t$-differentially stable:
$$
\sK=\Bbb C\biggl(t,z,\frac{\d z}{\d t},\dots,\frac{\d^6z}{\d t^6},
t_2,\frac{\d t_2}{\d t},\dots,\frac{\d^6t_2}{\d t^6},
t_3,\frac{\d t_3}{\d t},\dots,\frac{\d^6t_3}{\d t^6}\biggr)
$$
(provided that we are working locally,
in a neighbourhood of non-singular point $z=z_0$).
Finally, by~\thetag{23} we obtain both $F\in\sK$ and
$t_2,t_3\in\Bbb C[t,F,\d F/\d t]$, that is,
$$
\sK
=\Bbb C\biggl(t,z,\frac{\d z}{\d t},\dots,\frac{\d^6z}{\d t^6}\biggr)
\<F\>_{\d/\d t}
=\Bbb C(t)\<z,F\>_{\d/\d t},
$$
where $\sC\<F\>_{\d/\d t}$ is
the algebraic closure of the field
$$
\sC\biggl(F,\frac{\d F}{\d t},\frac{\d^2F}{\d t^2},
\frac{\d^3F}{\d t^3},\dots\biggr)
$$
(and $\sC$ is a $\d/\d t$-differentially stable field).
By~\thetag{24} and Proposition~3, the function $F(t)$ satisfies
some tenth order algebraic differential equation.

\proclaim{Proposition 6}
The transcendence degree over $\Bbb C(t)$ of the field
$$
\sK=\Bbb C(t)\<z,F\>_{\d/\d t}
$$
is~$10$.
\endproclaim

\demo{Proof}
By Proposition~2, the fields $\sL$ and $\sK$ coincide up to algebraic
extension. Thus,
$$
\trdeg_{\Bbb C(t)}\sK=\trdeg_{\Bbb C}\sK-1
=\trdeg_{\Bbb C}\sL-1=\trdeg_{\Bbb C(z)}\sL,
$$
which is the dimension of the differential Galois group of~\thetag{20}.
By~\cite{BH}, the Zariski closure of the projective monodromy
group of differential equation~\thetag{20}
(which is precisely its differential
Galois group) is $\SP_4(\Bbb C)$,
and $\dim_{\Bbb C}\SP_4$ is~$10$.
This completes the proof.
\enddemo

\remark{Remark}
Of course, we must scale $z\mapsto5^{-5}z$ to use the
results of~\cite{BH} above and later.
\endremark

\widestnumber\item{\kern1mm}
\proclaim{Theorem 2}
The differentially closed field
$$
\sK_1=\Bbb C(t)\<F\>_{\d/\d t}
$$
coincides \rom(up to algebraic extension\rom)
with the Picard--Vessiot extension~$\sL$
corresponding to the linear differential equation~\thetag{20}.
In particular\rom,
\roster
\item"(a)" $\trdeg_{\Bbb C(t)}\sK_1=10$;
\item"(b)" the mirror map $z(q)=z(e^t)$ is algebraic over~$\sK_1$.
\endroster
\endproclaim

\demo{Proof}
The proof is an immediate consequence of Propositions~3,~5, and~6.
\enddemo

We need the following functional version of the Schanuel conjecture
proved in~\cite{A}.

\proclaim{The Ax theorem}
Let $h_1(\tau),\dots,h_m(\tau)$ be analytic functions in
some neighbourhood of $\tau=0$ such that $h_1(\tau)-h_1(0)$
are linearly independent over~$\Bbb Q$. Then
$$
\trdeg_{\Bbb C}\Bbb C\bigl(h_1(\tau),\dots,h_m(\tau),
e^{h_1(\tau)},\dots,e^{h_m(\tau)}\bigr)\ge m+1,
$$
provided that all exponentials are defined.
\endproclaim

\proclaim{Proposition 7}
The function $q(t)=e^t$ is transcendent over~$\sK$.
\endproclaim

\demo{Proof}
To describe the projective monodromy group of~\thetag{20},
it is reasonable to scale
$$
t\mapsto\tau=\tau_1=\frac t{2\pi i},
\quad
t_2\mapsto\tau_2=\frac{t_2}{(2\pi i)^2},
\quad
t_3\mapsto\tau_3=\frac{t_3}{(2\pi i)^3},
\qquad
z\mapsto\wt z=5^{-5}z.
$$
Our aim is to prove the transcendence of $e^{\kappa\tau}$
over~$\sK$ for any $\kappa\in\Bbb C\setminus\{0\}$.

The simplest way to express the projective monodromy group
$G\subset\GL_4(\Bbb C)$ for the equation~\thetag{20}
in terms of~$\wt z$ is given by the Levelt theorem (see~\cite{BH}).
But in our case, studied in detail by physicists,
we can use the ready result from~\cite{H} for the local
monodromy $\gamma_1\in G$ about $\wt z=1$
(we adopt our basis near $\wt z=0$ with the same one from~\cite{H}):
$$
\pmatrix \tau_3 \\ \tau_2 \\ \tau_1 \\ 1 \endpmatrix
\mapsto
\gamma_1\pmatrix \tau_3 \\ \tau_2 \\ \tau_1 \\ 1 \endpmatrix,
\qquad
\gamma_1=\pmatrix 1 & 0 & 0 & 0 \\ 0 & 1 & 0 & 0 \\
0 & 0 & 1 & 0 \\ -1 & 0 & 0 & 1 \endpmatrix.
$$
In particular, we have
$$
\gamma_1^j:\tau\mapsto\frac\tau{1-j\tau_3}
$$
for any $j\in\Bbb Z$. We note that $\tau_3(\tau)$ is transcendent
over~$\Bbb C(\tau)$ by~\thetag{23} and Theorem~2.

Functions $\tau/(1-j\tau_3)$, $j=0,1,\dots,m$,
where $m\in\Bbb N$~is arbitrary,
are linearly independent over~$\Bbb C$
(and therefore, over~$\Bbb Q$).
By the Ax theorem, we get the algebraic independence over~$\Bbb C$
of at least~$m$ among $m+1$ exponentials
$$
\exp\biggl(\frac{\kappa\tau}{1-j\tau_3}\biggr), \qquad j=0,1,\dots,m.
\tag31
$$

We fix $m\ge12$. Assume now that $e^{\kappa\tau}$ is algebraic
over~$\sK$. Since the field~$\sK$ is $G$-in\-variant, acting by~$\gamma_1^j$
we see that all exponentials~\thetag{31} are algebraic over~$\sK$.
The transcendence degree over~$\Bbb C$ of the field generated
by exponentials~\thetag{31} is at least $m\ge12$, while
the transcendence degree over~$\Bbb C$ of~$\sK$ is~$11$ by Proposition~6.
This contradiction completes the proof of the transcendence
of~$e^{\kappa\tau}$ over~$\sK$.
\enddemo

\remark{Remark \rom1}
Our arguments using $\gamma_1\in G$ remain valid
for other linear differential equations producing mirror maps
(see~\cite{H}).
\endremark

\remark{Remark \rom2}
It is interesting to compare our approach of adding $q=e^{\kappa\tau}$
with the approach in~\cite{Ni}, where the author, roughly speaking,
considers the exponentials
$\exp(\kappa\tau/(1-j\tau))$ for $j=1,2,3,4$,
and proves their algebraic independence over $\Bbb C(\tau)$
(rather than~$\Bbb C$) from the
linear independence of $\tau/(1-j\tau)$, $j=1,2,3,4$, over~$\Bbb C$.
This proof is based on the Ostrowski--Kolchin theorem
(which is also mentioned in~\cite{A} as a related result).
\endremark

\remark{Remark \rom3}
D.~Bertrand (cf.~\cite{BZ2}) gives slightly different and easy
arguments proving Proposition~7
(and the Mahler--Nishioka result from~\cite{Ma},~\cite{Ni}).
Namely, when $\gamma$~runs through~$G$, the functions
$\exp(\kappa\cdot\gamma\tau)$ generate over~$\Bbb C$ a field
of infinite transcendence degree; therefore, they cannot be
algebraic over the field~$\sK$ of finite transcendence degree.
The key to the proof of this claim is the fact that
each of the functions $\exp(\kappa\tau/(1-j\tau_3))$
has essential singularities in the set
$S_j=\{\tau\in\Bbb C:\tau_3(\tau)=1/j\}$, $j\in\Bbb Z$,
thus it is transcendental over the field of functions
meromorphic in a neighbourhood of $s_j\in S_j$; in particular,
$\exp(\kappa\tau/(1-j\tau_3))$ is transcendental over the field
generated over~$\Bbb C$ by other exponentials
$\exp(\kappa\tau/(1-m\tau_3))$, $m\ne j$.
\endremark

Now, results on functional transcendence take their `close to final' form.

\proclaim{Theorem 3}
The transcendence degree over~$\Bbb C$ of the field
$\Bbb C(t,q=e^t)\<F\>_{\d/\d t}$ is~$12$.
\endproclaim

\proclaim{Theorem 4}
The transcendence degree over~$\Bbb C(q)$ of the field
$$
\sK_2=\Bbb C(q)\<K(q)\>_{\delta_q}
=\Bbb C(q)\<K(q)\>_{\d/\d q}
$$
is~$7$\rom, where $K(q)$ is the Yukawa coupling.
In other words\rom, $K(q)$ does not satisfy an algebraic
differential equation of order less than~$7$
with coefficients from $\Bbb C(q)$.
\endproclaim

\remark{Remark \rom1}
The same results can be received for other mirror
maps and Yukawa couplings considered in~\cite{M2} and~\cite{BS}.
\endremark

\widestnumber\item{\kern1mm}
\remark{Remark \rom2}
The result of Theorem~4 would be more useful if
there existed a set $\theta_1(q),\dots,\theta_7(q)$
of generators of~$\sK_2$ such that
\roster
\item"(i)" functions $\theta_1(q),\dots,\theta_7(q)$ have integral
coefficients of {\it polynomial\/} growth in their $q$-expansions,
\item"(ii)" these functions satisfy a (more or less) simple
system of non-linear differential equations.
\endroster
Condition~(i) is crucial for making number-theoretic applications,
while (ii) can simplify algebraic preliminaries.
\endremark

\remark{Remark \rom3}
It is natural to consider the differential field
$$
\sK_3=\Bbb C(q)\<z(q)\>_{\delta_q}
=\Bbb C(q)\<z(q)\>_{\d/\d q},
$$
since the seventh order differential equation for~$z(q)$ is simpler
than for~$K(q)$. However, in spite of the differential
duality observed in~\cite{KLRY} and~\cite{LY3},
we know of no arguments for~$K(q)$ to be algebraic over~$\sK_3$.
\endremark

\subhead
Acknowledgements
\endsubhead
This work is an extended version of author's talk
on the International Conference on Transcendental Numbers
held at Moscow, September 18--22, 2000.
I express my gratitude to all participants of this meeting
for their suggestive remarks.
Special gratitude is due to Professor D.~Bertrand
for his careful reading of the preliminary version of this paper
and for his valuable advice sharpening the description
of the transcendence results.
I am also grateful to Professor B.\,H.~Lian
for sending me the text of~\cite{LY3}.
This work was supported in part
by the INTAS--RFBR project no.~IR-97-1904.


\closegraphsfile
\enddocument
\end